\documentclass[12pt,reqno]{amsart} 
\usepackage{amssymb}
\usepackage{mathrsfs}
\usepackage{enumerate}
\usepackage[all]{xy}
\usepackage[usenames,dvipsnames]{color}
\usepackage[colorlinks=true, citecolor=OliveGreen, 
linkcolor=OliveGreen]{hyperref}
\usepackage{amsthm}

\renewcommand{\today}{\the\day/\the\month/\the\year}

\makeatletter
\@namedef{subjclassname@2020}{\textup{2020} Mathematics Subject 
Classification} 
\makeatother

\DeclareMathAlphabet\EuR{U}{eur}{m}{n}
\SetMathAlphabet\EuR{bold}{U}{eur}{b}{n}
\newcommand{\curs}{\EuR}
\newcommand{\catdef}[2][]{\expandafter\newcommand\csname#2\endcsname%
{#1\curs{#2}}}
\catdef{Ab}    \catdef{Grp}   \catdef{Cat}   \catdef{Vect}   \catdef{Mod}
\catdef{Is}    \catdef{Top}   \catdef{hoTop}  \catdef{ho}   \catdef{FibSimp}

\let\Gamma=\varGamma
\let\Omega=\varOmega
\let\Sigma=\varSigma

\renewenvironment{enumerate}[1][]
{\begin{enumerat}[#1]\setlength{\itemsep}{6pt}}{\end{enumerat}}

\newenvironment{enuma}{\begin{enumerate}[{\rm(a) }]}{\end{enumerate}}

\renewenvironment{itemize}
{\begin{itemiz}\setlength{\itemsep}{4pt} 
}{\end{itemiz}}

\definecolor{darkgreen}{rgb}{0,0.5,0}
\definecolor{bluegreen}{rgb}{0,0.2,0.8}
\definecolor{darkred}{rgb}{0.8,0,0}
\definecolor{newercolor}{rgb}{0.2,0,1}
\definecolor{darkyellow}{rgb}{0.7,0.7,0}
\definecolor{darkorange}{rgb}{0.8,0.4,0}

\newcommand{\mynote}[1]{{\color{blue}\noindent\textbf{\textup{[#1]}}}}

\numberwithin{table}{section}

\newlength{\short}
\newcommand{\4}[1]{\widebar{#1}}

\def\pair[#1,#2]{[\hskip-1.5pt[#1,#2]\hskip-1.5pt]}

\SelectTips{cm}{10} \UseTips   

\let\oldcirc=\circ
\renewcommand{\circ}{\mathchoice
    {\mathbin{\scriptstyle\oldcirc}}{\mathbin{\scriptstyle\oldcirc}}
    {\mathbin{\scriptscriptstyle\oldcirc}}
    {\mathbin{\scriptscriptstyle\oldcirc}}}

\def\beq#1\eeq{\begin{equation*}#1\end{equation*}}
\def\beqq#1\eeqq{\begin{equation}#1\end{equation}}

\numberwithin{equation}{section}

\newtheorem{Thm}{Theorem}[section]
\newtheorem{Prop}[Thm]{Proposition}

\newtheorem{Lem}[Thm]{Lemma}

\theoremstyle{definition}
\newtheorem{Defi}[Thm]{Definition}

\newcommand{\widebar}[1]
      {\overset{{\mskip1mu\leaders\hrule height0.4pt\hfill\mskip1mu}}{#1}
      \vphantom{#1}}

\newcounter{let} 
\loop\stepcounter{let}
\expandafter\edef\csname cal\alph{let}\endcsname%
{\noexpand\mathcal{\Alph{let}}}
\ifnum\thelet<26\repeat

\loop\stepcounter{let}
\expandafter\edef\csname scr\alph{let}\endcsname%
{\noexpand\mathscr{\Alph{let}}}
\ifnum\thelet<26\repeat

\newcommand{\tdef}[2][]{\expandafter\newcommand\csname#2\endcsname%
{#1\textup{#2}}}
\tdef{Iso}   \tdef{Aut}    \tdef{Out}    \tdef{Inn}    \tdef{Hom}
\tdef{End}   \tdef{Inj}    \tdef{map}    \tdef{Ker}    \tdef{Ob}
\tdef{Mor}   \tdef{Res}    \tdef{Id}     \tdef{Fr}     \tdef{Spin} 
\tdef{rk}    \tdef{conj}   \tdef{incl}   \tdef{proj}   \tdef{diag} 
\tdef{trf}   \tdef{Sol}    \tdef{He}     \tdef{Sz}     \tdef{cj}
\tdef{Rep}   \tdef{pr}    \tdef{Inndiag} \tdef{Outdiag}  \tdef{expt}
\tdef{supp}  \tdef{Isom}   \tdef{ord}    \tdef{Coker}   \tdef{Tr}
\tdef[_]{typ} \tdef[^]{op} \tdef[^]{ab}   \tdef{lcm}  \tdef{McL}
\tdef{restr}  \tdef{Comp}  \tdef{HS}     \tdef{ev}    \tdef{Stab}
\tdef{srk}

\newcommand{\fdef}[1]{\expandafter\newcommand\csname#1\endcsname%
{\mathfrak{#1}}}
\fdef{X}  \fdef{red}  \fdef{foc}  \fdef{hyp}  \fdef{Lie} \fdef{Y}

\newcommand{\bbdef}[1]{\expandafter\newcommand%
\csname#1\endcsname{\mathbb{#1}}}
\bbdef{C} \bbdef{F} \bbdef{R} \bbdef{Z} \bbdef{N} \bbdef{Q} \bbdef{K}

\newcommand{\itdef}[1]{\expandafter\newcommand\csname#1\endcsname%
{\textit{#1}}}
\itdef{PSL}  \itdef{PSU}  \itdef{SL}  \itdef{SU}  \itdef{GL} \itdef{GU}
\itdef{Sp}   \itdef{PSp} \itdef{PSO} \itdef{SO}   \itdef{SD} \itdef{PGU} 
\itdef{PGL}  \itdef{Co}  \itdef{Fi}  \itdef{GO}   \itdef{BDI}

\newcommand{\gee}{\varepsilon}

\newcommand{\lie}[3]{\def\test{#2}\def\tst{G}\ifx\test\tst{{}^{#1}#2_{#3}}
\else{{}^{#1}\!#2_{#3}}\fi}
\renewcommand{\*}{\,\lower6pt\hbox{\Large{\textup{*}}}\,}
\newcommand{\syl}[3][]{\textup{Syl}^{#1}_{#2}(#3)}
\newcommand{\sylp}[2][]{\syl[#1]{p}{#2}}

\newcommand{\defeq}{\overset{\textup{def}}{=}}

\newcommand{\mxfoura}[8]{\left(\begin{smallmatrix}#1&#2&#3&#4\\#5&#6&#7&#8}
\newcommand{\mxfourb}[8]{\\#1&#2&#3&#4\\#5&#6&#7&#8\end{smallmatrix}\right)}

\newcommand{\nsg}{\trianglelefteq}

\newcommand{\til}[1]{\widetilde{#1}}

\newcommand{\longleft}[1]{\;{\leftarrow%
\count255=0 \loop \mathrel{\mkern-6mu}%
    \relbar\advance\count255 by1\ifnum\count255<#1\repeat}\;}
\newcommand{\longright}[1]{\;{\count255=0 \loop \relbar\mathrel{\mkern-6mu}%
    \advance\count255 by1\ifnum\count255<#1\repeat\rightarrow}\;}

\newcommand{\RIGHT}[3]{\mathrel{\mathop{\kern0pt\longright#1}
        \limits^{#2}_{#3}}}

\newcommand{\LEFT}[3]{\mathrel{\mathop{\kern0pt\longleft#1}\limits^{#2}_{#3}}
}
\newcommand{\dRIGHT}[3]{\mathrel{%
   \mathop{\vcenter{\baselineskip=0pt\hbox{$\kern0pt\longright#1$}%
   \hbox{$\kern0pt\longright#1$}}}\limits^{#2}_{#3}}}
\newcommand{\LRIGHT}[3]{\mathrel{%
   \mathop{\vcenter{\baselineskip=0pt\hbox{$\kern0pt\longleft#1$}%
   \hbox{$\kern0pt\longright#1$}}}\limits^{#2}_{#3}}}
\newcommand{\RLEFT}[3]{\mathrel{%
   \mathop{\vcenter{\baselineskip=0pt\hbox{$\kern0pt\longright#1$}%
   \hbox{$\kern0pt\longleft#1$}}}\limits^{#2}_{#3}}}
\newcommand{\onto}[1]{\;{\count255=0 \loop \relbar\mathrel{\mkern-6mu}%
    \advance\count255 by1
    \ifnum\count255<#1 \repeat \twoheadrightarrow}\;}

\newcommand{\SFL}[1][]{(S#1,\calf#1,\call#1)}

\newcommand{\longline}{\bigskip\hfill\hbox to 8cm{\hrulefill}%
\hfill\bigskip}

\def\LFS(#1){\textup{LFS($#1$)}} 
\def\LF(#1){\textup{LF($#1$)}} 

\fdef{Fin}
\fdef{I}
\tdef{Ext}

\newcommand{\higherlim}[2]{\displaystyle\setbox1=\hbox{\rm lim}
	\setbox3=\hbox{$\scriptstyle{#1}$}
	\ifdim\wd1<\wd3
	\mathop{\hphantom{^{#2}}\vtop{\baselineskip=5pt\box1}^{#2}}_{#1}
	\else
	\mathop{\vtop{\baselineskip=5pt\box1}^{#2}}_{#1}
	\fi}

\newenvironment{ChgProp}[2]{\setcounter{section}{#1} \setcounter{Thm}{#2} 
\addtocounter{Thm}{-1} \begin{Prop}}{\end{Prop}}

\begin{document}

\title{Correction to: Finite approximations of $p$-local compact groups}

\author{Alex Gonzalez}
\address{Institut Arraona, Carrer de Praga, 43, 08207 Sabadell, Spain}
\email{agondem@gmail.com}
\thanks{}

\author{Bob Oliver}
\address{Universit\'e Sorbonne Paris Nord, LAGA, UMR 7539 du CNRS, 
99, Av. J.-B. Cl\'ement, 93430 Villetaneuse, France}
\email{bobol@math.univ-paris13.fr}
\thanks{B. Oliver is partially supported by UMR 7539 of the CNRS}




\begin{abstract}
We point out some minor errors in a paper by the first author, and explain why 
they do not affect any of the main results in the paper.
\end{abstract}


\bigskip

\maketitle

We note here some minor errors and misstatements in \cite{Gonzalez}, and 
describe the changes that are needed to fix the proofs of its main results. 
We refer to \cite[Definitions 1.3 and 1.8]{Gonzalez}, or to the earlier 
paper \cite{BLO6}, for the definitions of saturated fusion systems, centric 
linking systems, and transporter systems over a discrete $p$-toral group. 
Telescopic transporter systems and finite retraction pairs are defined in 
\cite[Definitions 2.1 and 2.3]{Gonzalez}.

\subsection{Linking systems} These are never defined explicitly in 
\cite{Gonzalez} (only centric linking systems). As used in this paper 
(especially in Subsection 3B), a linking system associated to a fusion 
system $\calf$ over $S$ is a transporter system $\call$ associated to 
$\calf$ (see \cite[Definition 1.8]{Gonzalez}) with the additional property 
that $E(P)=\gee_P(C_S(P))$ for each subgroup $P\le S$ fully centralized in 
$\calf$. The condition that $\calf^{cr}\subseteq\Ob(\call)$, which is 
usually included in the definition of a linking system, is \emph{not} 
assumed here. In fact, one consequence of Propositions 3.16 and 3.17 is 
that the linking subsystems $\call_i$ of Definition 3.9 and Proposition 
3.15 satisfy this last condition for $i$ large enough. 

\smallskip

\subsection{Definition 1.15} In the definition of a \emph{fine} unstable 
Adams operation, we need to require that the degree $\zeta$ of the 
operation satisfy $\zeta\ne\pm1$ (not just 
$\zeta\ne1$). The condition that $\zeta\ne-1$ (if $p=2$) is important in 
Hypotheses 3.7 and 3.14, to guarantee that $S$ is the union of the 
$S_i\defeq C_S(\Psi_i)$ (where $\Psi_i=\Psi^{p^i}$).

\smallskip

\subsection{Lemmas 2.2 and 3.5, Remark 3.6, and Proposition 2.6} The 
following changes are needed:
\begin{enumerate}[$\bullet$] 

\item Drop Lemma 2.2.

\item In Lemma 3.5, Remark 3.6, and Proposition 4.3, add to each 
statement the extra assumption that the approximation $\{\SFL[_i]\}_{i\ge0}$ 
to $\calg=\SFL$ be taken with respect to a telescopic transporter system 
\emph{constructed using a finite retraction pair}.

\item In the proof of Lemma 3.5, replace the reference to Lemma 2.2 by one 
to Proposition 2.6.

\end{enumerate}
We do not know whether Lemma 2.2 holds after 
changing the conclusion to say that the inclusion of $|\call|$ in 
$|\til\call|$ is a mod $p$ homology equivalence (which would suffice to 
prove Lemma 3.5). The proof in \cite{Gonzalez} of this lemma is based on 
\cite[Proposition A.9]{BLO6}, whose hypotheses need not hold under the 
given assumptions.

\smallskip

\subsection{Proposition 2.7.} An additional assumption is needed in the 
statement of the proposition: that $\psi(P)^\star=\psi(P^\star)$ for every 
$P$ in $\Ob(\calt)$ (where $\psi\in\Aut(S)$ is the automorphism induced by 
$\Psi\in\Aut\typ^I(\call)$). With this extra condition, the proof holds 
without modification. This proposition is used in \cite{Gonzalez} in the 
following two situations: 
\begin{enumerate}[$\bullet$]

\item If $\calt$ is a centric linking system, $\Psi$ is an isotypical 
automorphism of $\call$, and $\til{\calt}$ is the telescopic linking system 
defined via the ``bullet construction'' as in Example 2.8 (i.e., 
$(-)^\star=(-)^\bullet$), then the condition 
$\psi(P)^\bullet=\psi(P^\bullet)$ always holds. In particular, this applies 
to all unstable Adams operations defined on $\call$ used in the paper.

\item In Step 3 of the proof of Proposition 3.17, one needs to consider the 
case where $\calt$ is a linking system, $A \leq T$ is a subtorus, $\calt = 
N_\call(A)$ is the normalizer linking system of $A$, and $\til{\calt}$ is 
the telescopic linking system associated to $\call$ via the retraction pair 
$((-)^\star_{N_\calf(A)},(-)^\star_{N_\call(A)})$ of Example 2.9. 
Proposition 2.7 is applied with $\Psi$ the restriction of an automorphism 
of $\call$ (an unstable Adams operation), and where $P^\star=P^\bullet$ for 
all $P\in\Ob(N_\call(A))$. The extra assumption 
$\psi(P)^\star=\psi(P^\star)$ then holds, and Proposition 2.7 can still be 
applied.

\end{enumerate}

\noindent This also means that in Hypothesis 3.7 (hence throughout 
Subsection 3A), we must assume that $(-)^\star_\calf=(-)^\bullet_\calf$ and 
$(-)^\star_\call=(-)^\bullet_\call$ (or at least that 
$\Psi(P^\star)=\Psi(P)^\star$ for all $P\in\Ob(\call)$).

\smallskip

\subsection{Definition 3.1 and Lemma 3.8} In Definition 3.1, we need to 
include the additional assumption that the sequences $\{S_i\}_{i\ge0}$, 
$\{P_i\}_{i\ge0}$, and $\{Q_i\}_{i\ge0}$ are all increasing. These extra 
conditions do hold when finite approximations are constructed in the proof 
of Theorem 1: $S_i\le S_{i+1}$ for each $i$ since $\Psi_{i+1}=(\Psi_i)^p$ 
by Hypothesis 3.7, while $P_i=P\cap S_i$ and $Q_i=Q\cap S_i$ for each $i$.

Likewise, in Lemma 3.8, the conclusion that the sequence $\{S_i\}$ is 
increasing should be added to point (i). This is needed in the proofs of 
Proposition 3.11(i,ii), Proposition 3.17 (Step 6), and Theorem 1.

\smallskip

\subsection{Proposition 3.15} The claim that each object in $\call_i$ be 
$\calf_i$-centric should be dropped from this statement. In other words, 
the proposition should only state that $\call_i$ is a linking system for 
all $i$. In the fourth line of the proof, the last equality 
(``$\,=\gee_i(Z(P))$'') should be removed. Proposition 3.15 is used only in the proof of Theorem 1, and only in this 
restricted form.

\smallskip

\subsection{Proposition 3.16} The statement of this proposition should be 
modified as follows:

\begin{ChgProp}3{16}
Let $i\ge0$, and let $P\le S_i$ be such that $P\lneqq P^\bullet\cap S_i$. 
Then either 
\begin{enuma} 

\item $P$ is not $\calf_i$-centric; or 

\item $P$ is $\calf_i$-centric and 
	$ \Out_{S_i}(P) \cap O_p(\Out_{\calf_i}(P)) \ne 1. $

\end{enuma}
In particular, $P$ cannot be both $\calf_i$-centric and $\calf_i$-radical.
\end{ChgProp}

\noindent Thus the inequality $\Out_{S_i}(P)\cap 
O_p(\Out_{\calf_i}(P))\ne1$ is shown only when $P$ is $\calf_i$-centric. 
The proof of the proposition remains unchanged. Its only application, in 
the proof of Theorem 1, is unaffected.

\smallskip

\subsection{Proof of Proposition 3.17} Throughout this proof, $T$ always denotes the maximal torus of the Sylow group $S$. Also, for each $P\le S$, $T_P\nsg P$ denotes its maximal torus.

There is a typo in the statement of (\ddag): it should end with 
``some $Q\le S_i$ such that $K_Q\cap\Aut_{S_i}(Q)$ contains some element 
which is not in $\Inn(Q)$'' (\emph{not} ``such that $K_Q\cap\Aut_S(Q)$ 
\dots''). However, this does not affect the rest of the proof. The precise 
formulation of (\ddag) is used only in Step 3 (statement (3-h)) and at 
the end of Steps 6 and 8. Most of the other steps are of the form ``if 
(\ddag) holds for subgroups of one type, then it holds for those of another 
type'', and the difference between the two versions of (\ddag) doesn't 
affect those arguments.

In the next-to-last paragraph in Step 10, there is a gap in the proof that 
$K'_j\cap S_i\in\Ob(\call_i)$ for each $1\le j\le n$. One way to show this 
is as follows. Point (10-b), 
together with the same argument applied to $\4{x}^{-1}$, shows that 
$\4{K}'_j\cap\4{S}_i$ is $\4{T}$-conjugate to $\4{K}_j\cap\4{S}_i$. So by 
(4-a), this implies that $K'_j\cap S_iA$ is $T$-conjugate to $K_j\cap 
S_iA$. Then by the lemma below, 
	\[ (K'_j\cap S_i)^\bullet=(K'_j\cap S_iA)^\bullet 
	\qquad\textup{and}\qquad 
	(K_j\cap S_i)^\bullet=(K_j\cap S_iA)^\bullet, \]
and hence the two are $T$-conjugate. So $(K'_j\cap 
S_i)^\bullet\in\Ob(\call)$ since $(K_j\cap S_i)^\bullet\in\Ob(\call)$. 

\noindent\textbf{Lemma.} \emph{Assume $P$ is such that $A\le P\le S$. If $p=2$, 
assume $\Psi_i$ has degree $1$ mod $4$. Then $(P\cap S_iA)^\bullet = (P\cap 
S_i)^\bullet$.} 

\begin{proof} In the notation of Definition 1.12, we have 
	\[ (P\cap S_i)^{[e]}\cdot A^{[e]} \le (P\cap S_iA)^{[e]} 
	\le (P\cap S_i)^{[e]}\cdot A, \]
with equality since $A=A^{[e]}$ ($A$ is infinitely divisible). Set 
$W=\Aut_\calf(T)$. Set $\Lambda=C_W((P\cap S_i)^{[e]})$, and 
$\Lambda_0=C_W(P\cap S_iA)^{[e]} = \Lambda\cap C_W(A)$. Then 
$\Lambda=\Lambda_0$ since $C_W(A\cap S_i)=C_W(A)$ by the assumption on the 
degree of $\Psi_i$. So 
	\[ A \le I((P\cap S_iA)^{[e]}) = I((P\cap S_i)^{[e]}), \]
and hence 
	\begin{equation*} 
	(P\cap S_iA)^\bullet = (P\cap S_i)\cdot A\cdot I((P\cap S_iA)^{[e]})_0
	= (P\cap S_i)\cdot I((P\cap S_i)^{[e]})_0 = (P\cap S_i)^\bullet. 
	\qedhere
	\end{equation*}
\end{proof}

\bigskip

\smallskip

\subsection{Proof of Lemma 5.3.} Step 2 can be simplified by noting that 
$\Ker(\omega_P)$ is discrete $p$-toral by \cite[Lemma 1.7]{Gonzalez}. Since 
$\Ker(\rho_P)=\gee_P(C_S(P))$ is also discrete $p$-toral, so is 
$\Ker(\4\rho_{P/V})=E(P/V)$. Also, $P/V$ is fully centralized since it is 
$\calf/V$-centric, so we have $\gee_{P/V}(Z(P/V))\in\sylp{E(P/V)}$, 
finishing the proof that $E(P/V)=\gee_{P/V}(Z(P/V))$ and hence 
$(\til\call/V)^c$ is a linking system. With this simplification, Step 1 can 
be removed. 

In Step 3, to apply \cite[II.5.1]{BK}, one must first show that the action 
of $\pi_1(|(\til\call/V)^c|)$ on $H_*(BV;\F_p)$ is nilpotent. 
Alternatively, one can instead apply the principal fibration lemma 
\cite[Lemma II.2.2]{BK}, after showing that the fibration sequence is 
equivalent to that of a principal fibration.


\end{document}